\documentclass[10pt]{amsart}
\usepackage{latexsym}
\usepackage{amscd, amsfonts, eucal, mathrsfs, amsmath, amssymb, amsthm}
\usepackage[OT2,OT1]{fontenc}

\newcommand{\field}[1]{\mathbf #1}
\newcommand{\mf}[1]{\mathfrak #1}
\newcommand{\mc}[1]{\mathcal #1}
\newcommand{\ms}[1]{\mathscr #1}
\newcommand{\widebar}[1]{\overline{#1}}

\newcommand{\Z}{\field Z}

\newcommand{\simto}{\stackrel{\sim}{\to}}

\renewcommand{\phi}{\varphi}

\newcommand{\br}{\mf{Br}}

\renewcommand{\hom}{\operatorname{Hom}}
\newcommand{\shom}{\ms H\!om}
\DeclareMathOperator{\chom}{\underline{Hom}}

\newcommand{\spec}{\operatorname{Spec}}

\newcommand{\rspec}{\operatorname{\bf Spec}}

\renewcommand{\P}{\field P}
\newcommand{\A}{\field A}

\DeclareMathOperator{\uquot}{\underline{Quot}}
\DeclareMathOperator{\quot}{Quot}

\newcommand{\GL}{\operatorname{GL}}

\DeclareMathOperator{\ext}{\operatorname{Ext}}

\DeclareMathOperator{\coh}{\operatorname{Coh}}
\DeclareMathOperator{\vbdl}{\operatorname{LF}}

\DeclareMathOperator*{\tensor}{\otimes}

\DeclareMathOperator{\Tr}{\operatorname{Tr}}

\newcommand{\surj}{\twoheadrightarrow}
\newcommand{\inj}{\hookrightarrow}
\newcommand{\id}{\operatorname{id}}

\DeclareMathOperator{\coker}{\operatorname{coker}}

\DeclareMathOperator{\isom}{\operatorname{Isom}}

\DeclareMathOperator{\B}{\operatorname{B\!}}

\newtheorem{lem}{Lemma}[section]
\newtheorem{thm}[lem]{Theorem}

\newtheorem{cor}[lem]{Corollary}

\newtheorem{lems}{Lemma}[subsection]
\newtheorem{thms}[lems]{Theorem}
\newtheorem{props}[lems]{Proposition}

\theoremstyle{definition}
\newtheorem{defn}[lem]{Definition}

\newtheorem{example}[lem]{Example}

\newtheorem{paras}[lems]{}

\theoremstyle{remark}
\newtheorem{remark}[lem]{Remark}
\newtheorem{remarks}[lem]{Remarks}

\newtheorem{remarksub}[lems]{Remark}

\title{Remarks on the stack of coherent algebras}
\author{Max Lieblich}

\begin{document}

\maketitle

\begin{abstract}
  We consider the stack of coherent algebras with proper support, a
  moduli problem generalizing Alexeev and Knutson's stack of
  branchvarieties to the case of an Artin stack.  The main results are
  proofs of the existence of Quot and Hom spaces in greater generality
  than is currently known and several applications to Alexeev and
  Knutson's original construction: a proof that the stack of
  branchvarieties is always algebraic, that limits of one-dimensional
  families always exist, and that the connected components of the
  stack of branchvarieties are proper over the base under certain
  hypotheses on the ambient stack.
\end{abstract}
\section{Introduction}

In a recent paper \cite{alexeev-knutson}, Alexeev and Knutson consider
a moduli problem closely related to the Hilbert scheme: the stack of
branchvarieties.  They focus on branchvarieties over a projective
scheme, and they prove (among other things) that the stack has proper
components.  The question of the existence of limits along discrete
valuation rings is left open for not necessarily projective schemes,
or, more generally, for algebraic spaces, as is the question of
quasi-compactness of the connected components of the stack.

A branchvariety over $Y$ is, in particular, a finite morphism.
Abandoning the projective methods of \cite{alexeev-knutson}, we give
an approach to the stack of branchvarieties using standard abstract
methods applied to the stack of all finite morphisms with proper
support.  In particular, we give an abstract ``construction'' of (a
mild generalization of) Alexeev and Knutson's stack of branchvarieties
over any Deligne-Mumford stack.  Then we show that the stack satisfies
the valuative criterion of properness in full generality.  When the
stack has quasi-projective coarse moduli space, we show that the
connected components of the stack of branchvarieties are proper.  In
fact, one can give (somewhat contrived) explicit numerical constraints
which ensure boundedness of a collection of branchvarieties, if
desired.

Among other things, our methods also yield proofs of the algebraicity
of the usual stacks under very general conditions: $\quot$ spaces for
quasi-coherent sheaves of finite presentation with proper support on
an Artin stack of finite presentation and (as pointed out to us by
Starr) $\hom$ spaces from a proper Artin stack of finite presentation
to a global quotient by a flat linear algebraic group scheme.  This
gives generalizations of various recent results of Starr and Olsson
contained in \cite{olsson-homs}, \cite{olsson-starr}, and
\cite{starr}.  In fact, \ref{P:quotient-homs} gives a natural
complement to Theorem 1.1 of \cite{olsson-homs} for target stacks with
potentially positive-dimensional stabilizers (but which are required
to be global quotients).  These applications are taken up in
\S{~\ref{S:construction}} and its subsections.

As the initial motivation for this work was to understand Alexeev and
Knutson's stack, we use the rest of this introduction to sketch
the main results concerning branchvarieties.  Throughout, $Y\to S$ is
a Deligne-Mumford stack of finite presentation over an excellent
algebraic space.  Somewhat more general results are obtainable in
certain cases (as in the recent work of Starr \cite{starr}), but we
restrict ourselves to this case because it is easier to
describe and has some potentially interesting applications.

\begin{defn}
  A \emph{family of branchvarieties\/} over $Y$ parametrized by an
  affine scheme $T\to S$ is a proper flat Deligne-Mumford stack of
  finite presentation $X\to T$ with geometrically reduced fibers
  equipped with a finite $T$-morphism $X\to Y_T$.
\end{defn}

Note that when $Y$ has automorphisms, we allow $X$ to have
automorphisms.  The locus consisting of branchvarieties with $X$ an
algebraic space defines an open substack, but properness requires that
we allow $X$ to acquire automorphisms in the limit.  (Indeed, one can
easily make examples where the presence of automorphisms is necessary
in the limit by considering the Hilbert scheme of a weighted
projective stack.)  Of course, when $Y$ is an algebraic space (resp.\
scheme), $X$ will also be an algebraic space (resp.\ scheme).

We will prove the following.

\begin{thm}\label{T:existence}
  Branchvarieties over $Y$ form an Artin stack $\br_{Y/S}$ locally of
  finite presentation over $S$.
\end{thm}

The most natural way to build the stack is as a stack of flat families
of algebras over the structure sheaf of $Y$, something we sketch in
section \ref{S:construction}.

\begin{thm}\label{T:limits}
  Let $T\to S$ be a discrete valuation scheme with generic point
  $\eta$.  Suppose $X_{\eta}\to Y_{\eta}$ is a branchvariety.  There
  is a finite totally ramified extension $T'/T$ and a family of
  branchvarieties $X\to Y_{T'}$ extending the generic family.  Given
  two such families $X_1$ and $X_2$ and an isomorphism
  $\phi:(X_1)_{\eta}\to (X_2)_{\eta}$ over $Y_{\eta}$, there is a
  unique isomorphism $X_1\to X_2$ over $Y$ extending $\phi$.
\end{thm}

The key to proving \ref{T:limits} is to note that both parts -- the
existence and uniqueness -- are local in the \'etale topology.
Uniqueness then allows us to glue the local solutions to get a global
limit.

\begin{thm}\label{T:properness}
  If $S$ is Noetherian and either 1) $Y$ admits a proper flat cover by
  a quasi-projective algebraic space over $S$, or 2) $Y$ is a quotient
  stack whose coarse moduli space is quasi-projective over $S$, then
  $\br_{Y/S}$ is a disjoint union of proper $S$-stacks.
\end{thm}

Having proven that limits exist, the point is to prove that the
connected components are quasi-compact.  This is proved using the
various numerical invariants provided by polarizations of coverings or
of the coarse moduli space.

\section*{Acknowledgments}
\label{sec:acknowledgments}

I would like to thank Jason Starr and the referee for many helpful
comments and suggestions, and in particular for strongly suggesting
inclusion of sections \ref{sec:digr-quot-spac} and
\ref{sec:digr-hom-stacks}.

\section{The stack of coherent algebras with proper support}
\label{S:construction}

In this section, we sketch one approach to studying the stack of
finite morphisms (and ultimately proving \ref{T:existence}).  
The recent preprint of Starr \cite{starr} provides a general method
for considering a whole slew of relatively algebraic morphisms of
stacks.  However, as Starr has pointed out to us, in the case of
branchvarieties his method will only produce the open substack of our
construction consisting of branchvarieties which are algebraic spaces.
It is thus with the goal of producing a full compactification of the
moduli problem that we present the following method.

\subsection{Algebraicity of the stack}
\label{S:alg-funds}

In this section, we show that the collection of finite
morphisms with proper support to an Artin stack of finite presentation
over an excellent algebraic space is itself an Artin stack locally of
finite presentation over the base.

\begin{thms}\label{T:coh-rep}
  Let $Y\to S$ be an Artin stack of finite presentation over
  an excellent algebraic space.  The stack $\coh_{Y/S}$ of flat
  families of finitely presented quasi-coherent sheaves on $Y$ with
  proper support over $S$ is an Artin stack locally of finite
  presentation over $S$.
\end{thms}
\begin{proof}
  This is a straightforward application of Artin's theorem, using the
  stacky Grothendieck existence theorem proven by Olsson in 
  \cite{olsson} and the deformation theory of Illusie for
  modules in a topos.  When $Y$ is a Deligne-Mumford stack, one can
  simply use the \'etale topos of $Y$.  However, when $Y$ is an
  arbitrary Artin stack, one must work with Cartesian modules on a
  simplicial topos generated by a smooth cover of $Y$.  (To apply the
  theory of Illusie one must in addition note that any infinitesimal
  deformation of a Cartesian module is Cartesian; this ensures that
  the full $\ext$ groups parametrize obstructions and deformations.)
  The reader unfamiliar with these methods is referred to
  \cite{olsson-sheaves} for further details (but should note that
  deformation theory of sheaves is not explicitly developed there).
  The uncomfortable reader may choose to only think about
  Deligne-Mumford stacks.
\end{proof}

\begin{remarksub}
  This gives a generalization of Th\'eor\`eme 4.6.2.1 of \cite{lmb}, which only
  handles the case of projective morphisms $Y\to S$ of Noetherian schemes.
\end{remarksub}

\begin{props}\label{P:hom-rep}
  Let $Y\to S$ be a morphism of finite presentation between Artin
  stacks.  Let $\ms F$ and $\ms G$ be finitely presented
  quasi-coherent sheaves on $Y$ such that $\ms G$ is $S$-flat and the
  support of $\ms G$ is proper over $S$.  The stack of homomorphisms
  $\chom(\ms F,\ms G)\to S$ is representable by algebraic spaces
  locally of finite type.
\end{props}
\begin{proof}
  The stack $\chom(\ms F,\ms G)$ is the stack on the big fppf site of
  $S$ whose fiber category over a $1$-morphism from a scheme $W\to S$
  is the set of homomorphisms $\hom(\ms F_W,\ms G_W)$.  It is easy to
  check that this is a sheaf on the big fppf topos of $S$; to show
  that it is a relative algebraic space it suffices to prove this
  after pulling back to a smooth cover of $S$, so we may assume that
  $S$ is a scheme.  Now we can apply Artin's theorem (for algebraic
  spaces!), but using elementary arguments in place of Illusie's
  theory.  E.g., given a homomorphism from $\ms F_W\to\ms G_W$ and an
  infinitesimal extension of $W$ by an ideal $I$, basic homological
  algebra yields an obstruction to lifting the homomorphism in
  $\ext^1(\ms F,I\tensor\ms G)$.  The set of lifts is a torsor under
  $\hom(\ms F,I\tensor\ms G)$ (again by basic algebra).  The
  $S$-flatness of $\ms G$ is used in verifying the Schlessinger
  conditions on the functor $\chom(\ms F,\ms G)$.  The algebraic
  approximation of formal deformations is accomplished using the
  Grothendieck Existence Theorem for Artin stacks proven by Olsson
  \cite{olsson}. The rest is completely straightforward.
\end{proof}

\begin{remarksub}
  One easily deduces from Corollaire 7.7.8 of \cite{ega3-2} that
  \ref{P:hom-rep} holds when $Y\to S$ is a projective morphism of
  locally Noetherian schemes, and moreover that $\chom(\ms F,\ms G)$
  is representable locally on $S$ by the kernel of a linear
  homomorphism of geometric vector bundles.  The method above already
  gives a generalization when $Y\to S$ is a proper morphism of finite
  presentation between schemes, something which seems difficult to
  prove with the techniques of \cite{ega3-2}, except when $\ms F$ is
  the cokernel of a homomorphism between locally free sheaves of
  finite rank.
\end{remarksub}

\begin{props}
  Let $Y\to S$ be a morphism of finite presentation of Artin stacks
  and $\ms F$ an $S$-flat finitely presented quasi-coherent sheaf on
  $Y$ with proper support.  The stack of commutative algebra
  structures on $\ms F$ is represented by a relative algebraic space
  of finite presentation over $S$.
\end{props}
\begin{proof}
  To give an algebra structure on $\ms F$ is to give 1) a map $\mu:\ms
  F\tensor\ms F\to\ms F$ and 2) a map $\upsilon:\ms O\to\ms F$ such
  that $\mu$ defines a commutative and associative pairing and
  $\upsilon$ defines a unit for this pairing.  The conditions of
  commutativity, associativity, and unit define a closed substack of
  the stack of pairs $\chom(\ms F\tensor\ms F,\ms F)\times\chom(\ms
  O,\ms F)$.  Since $\ms F$ is $S$-flat, the latter stack is
  algebraic.  The result follows.
\end{proof}

\subsection{$\quot$ spaces}
\label{sec:digr-quot-spac}

Using the above methods, one can prove that algebraic $\quot$ spaces
exist on Artin stacks, generalizing the main result of \cite{olsson-starr}.

\begin{paras}
  Let $Y\to S$ be a morphism of finite presentation between Artin
  stacks and $\ms F$ and $\ms G$ quasi-coherent sheaves of finite
  presentation on $Y$ such that $\ms G$ is $S$-flat and has proper support.
  We may apply \ref{P:hom-rep} to produce an Artin stack with a
  representable morphism $\chom(\ms F,\ms G)\to S$ locally of finite
  presentation which parametrizes homomorphisms $\ms F\to\ms G$.
\end{paras}
\begin{lems}\label{L:surj-part}
  There is an open substack of $\chom(\ms F,\ms G)$ parametrizing surjective
  homomorphisms.
\end{lems}
\begin{proof}
  By Nakayama's lemma (and the faithful flatness of field extensions),
  given a scheme $W\to S$, a homomorphism $\phi:\ms F_W\to\ms G_W$ is
  surjective if and only if it is surjective when pulled back to all
  geometric points of $W$.  Thus, to show that surjections are
  represented by an open substack, it suffices to show that when $W$ is
  reduced the locus of points in $W$ over which $\phi$
  is surjective is open.  Since the formation of $\coker\phi$ commutes
  with base change, this is the same as showing that given a
  quasi-coherent sheaf $\ms K$ on $Y_W$ of finite presentation whose
  support is proper over $W$, the locus of points $w\in W$ such that
  $\ms K_w=0$ is open.  By Nakayama's lemma, the set of such $t$ is
  precisely the complement of the image of the support of $\ms K$.
  Since the support is proper, it has closed image.
\end{proof}

\begin{paras}   
  Let $Y\to S$ be a morphism of finite presentation between Artin
  stacks and $\ms F$ a quasi-coherent sheaf of $\ms O_Y$-modules of
  finite presentation.  We can define a presheaf on the category of
  $S$-schemes as follows: For any $S$-scheme $W\to S$, let
  $\quot_{Y/S}(\ms F)(W)$ be isomorphism classes of surjective
  homomorphisms $\ms F_W\surj\ms G$ with $\ms G$ a $W$-flat
  quasi-coherent sheaf on $Y_W$ of finite presentation with support
  proper over $W$.  Since sheaves and homomorphisms glue in the fpqc
  topology on an Artin stack, it is easy to see that $\quot_{Y/S}(\ms
  F)$ is a sheaf on the big fpqc topos on $S$.  In the parlance of
  \S{14} of \cite{lmb} (slightly generalized to base stacks rather
  than base schemes), the formation of $\quot_{Y/S}(\ms F)$ is a
  ``construction locale.''
\end{paras}
\begin{props}\label{P:quot-exists}
  There is a representable morphism locally of finite presentation
  $\uquot_{Y/S}(\ms F)\to S$ whose associated sheaf on the big fpqc
  site of $S$ is isomorphic to $\quot_{Y/S}(\ms F)$.
\end{props}
\begin{proof}
  By standard limiting methods, it suffices to prove the result when
  $S$ is an excellent algebraic space.  Let $\ms G$ be the universal
  sheaf on $Y\times\coh_{Y/S}$ and consider the $\coh_{Y/S}$-stack
  $\chom_{Y\times_S\coh_{Y/S}}(\ms F_{\coh_{Y/S}},\ms G)$, which is an
    Artin stack locally of finite presentation by
    \ref{T:coh-rep} and \ref{P:hom-rep}.  By \ref{L:surj-part} there
    is an open substack parametrizing surjective morphisms, and this
    is precisely the stack $\uquot_{Y/S}(\ms F)$.
\end{proof}

\subsection{$\hom$ stacks}
\label{sec:digr-hom-stacks}

Jason Starr has pointed out to us that we can refine the results of
\S{~\ref{S:alg-funds}} in another direction to study $\hom$-stacks
$\hom(X,Y)$ where $X$ is a proper Artin stack and $Y$ is a suitable
quotient stack.  


Let $S$ be an excellent algebraic space and $Y\to S$ a proper Artin
stack of finite presentation.

\begin{lems}\label{L:-loc-free}
  There is an open substack $\vbdl_{Y/S}\subset\coh_{Y/S}$
  parametrizing locally free sheaves of finite rank.
\end{lems}
\begin{proof}
  This follows immediately from Nakayama's lemma and reduction to the
  Noetherian case.
\end{proof}

Given an Artin $S$-stack $X\to S$, one can define a $\hom$-stack
$\hom_S(Y,X)$ to have objects over $W\to S$ the groupoid of
$1$-morphisms $Y_W\to X_W$.  Choosing a flat presentation $V\to Y$,
the objects of $\hom_S(Y,X)$ over $W$ are the same thing as
simplicial objects of $X_W$ over the simplicial scheme coming from the
collection fiber products $V_W\times_{Y_W}\cdots\times_{Y_W} V_W$.  (In other
words, passing to the associated simplicial object defines a natural
equivalence of groupoids.)  From this point of view, it is clear that
$\hom_S(Y,X)$ is an fppf $S$-stack.

One such stack we can take for $X$ is $\B{\GL_n}$, the classifying
stack of $\GL_n$-torsors.

\begin{lems}\label{L:hom-to-bgl}
  The stack $\hom_S(Y,\B{\GL_n})$ is an Artin stack locally of finite
  presentation over $S$.
\end{lems}
\begin{proof}
  It is easy to check that giving an object of $\hom_S(Y,\B{\GL_n})$ over
  $W\to S$ is the same as giving the $\GL_n$-torsor on $Y_W$
  associated to a locally free sheaf of $\ms O_{Y_W}$-modules of rank
  $n$.  Applying \ref{L:-loc-free} yields the result.
\end{proof}

\begin{lems}\label{L:sections}
  Let $X\to S$ and $Y\to S$ be Artin stacks locally of finite
  presentation over $S$ with $Y$ proper and $X$ separated and let
  $\phi:X\to Y$ be a representable $S$-morphism.  There is an Artin
  stack $\Sigma(\phi)$ locally of finite presentation over $S$ whose
  objects over $W\to S$ are pairs $(\psi,\gamma)$ with $\psi:Y_W\to
  X_W$ a $W$-morphism and $\gamma:\phi\psi\simto\id$ an isomorphism.
\end{lems}
\begin{proof}
  First, note that since $X\to Y$ is representable, it is easy to see
  that the stack $\Sigma(\phi)$ is the stack associated to a sheaf on
  $S$.

  We may assume that $S$ is Noetherian.  Consider the $\quot$ space
  $\quot_{X/S}(\ms O_X)$, whose points correspond to flat families of
  closed substacks of $X$.  The universal family $\mc Z\subset
  X\times_S\quot_{X/S}(\ms O)$ and the map $X\to Y$ give rise to a
  representable morphism $\mc Z\to Y\times_S\quot_{X/S}(\ms O)$ of
  proper $\quot_{X/S}(\ms O)$-stacks.  Taking the Stein factorization
  yields a finite coherent $\ms O_{Y_{\quot}}$-algebra $\ms A$.  The
  locus of $Y_{\quot}$ over which the natural map $\ms
  O_{Y_{\quot}}\to\ms A$ is an isomorphism is an open substack of
  $Y_{\quot}$.  The complement of its image in $\quot_{X/S}$ yields an
  open substack $\ms U\subset\quot_{X/S}(\ms O)$ over which the map
  $\mc Z\to Y_{\quot_{X/S}(\ms O)}$ is an isomorphism.  It is
  immediate that $\Sigma(\phi)$ is represented by $\ms U$.
\end{proof}

With this in hand, we can prove algebraicity of $\hom$-stacks to
global quotients by linear algebraic group schemes.

\begin{props}\label{P:quotient-homs}
  Let $X=[Z/G]$ be a quotient stack with $Z$ separated and of finite
  presentation over $S$ and $G$ an $S$-flat linear algebraic group
  scheme.  The stack $\hom_S(Y,X)$ is an Artin stack locally of finite
  presentation over $S$.
\end{props}
\begin{proof}
  Since algebraicity is \'etale-local, we may assume that $S$ is an
  excellent quasi-compact scheme.  There is an inclusion $G\inj\GL_n$
  for some $n$.  Since $G$ is $S$-flat, the quotient
  $Z'=Z\times^G\GL_n$ is an algebraic space of finite presentation
  over $S$ and there is a natural isomorphism $X=[Z'/\GL_n]$.  Thus,
  we may assume that $G=\GL_n$.  A map from $Y$ to $X$ is thus the
  same thing as an equivariant map from a $\GL_n$-torsor $T$ on $Y$
  to $Z'$, which is the same thing as a section of the bundle
  $Z'_Y\times^{\GL_{n,Y}}T\to Y$.

  There is a universal $\GL_n$-torsor $T\to
  Y\times_S\hom_S(Y,\B{\GL_n})$ which is representable by algebraic
  spaces.  Forming the $Z'$-bundle
  $Z'_{\hom_S(Y,\B{\GL_n})}\times^{\GL_n}T\to
  Y_{\hom_S(Y,\B{\GL_n})}$, we can apply \ref{L:sections} to see that
  the natural map $\hom_S(Y,X)\to\hom_S(Y,\B{\GL_n})$ is representable
  by algebraic spaces.
\end{proof}

\section{Applications to branchvarieties}
\label{sec:appl-branchv}

\subsection{The proof of \ref{T:existence}}
\label{sec:conclusion-proof}
Given the results of \S\ref{S:alg-funds}, the proof of \ref{T:existence}
follows immediately from the following lemma concerning the locus of
reduced fibers.

\begin{lems}
  Let $\pi:Y\to S$ be a proper flat morphism of finite presentation of
  Artin stacks.  The locus over which the geometric fibers of $Y$ are
  reduced is an open substack of $S$.
\end{lems}
\begin{proof} Let $\rho:U\to Y$ be a smooth cover.  By 12.2.1 of
  \cite{ega4-3}, the locus of points in $u\in U$ such that
  $U_{\pi(u)}$ is geometrically reduced at $u$ is open in $U$.  Taking
  the image of this open subscheme by $\rho$ defines an open substack
  $V\subset Y$.  Since $\pi$ is proper, $\pi(Y\setminus V)\subset S$
  is closed.  Since the property of being geometrically reduced is
  local in the smooth topology, it is easy to see that
  $S\setminus\pi(Y\setminus V)$ is the open set parametrizing
  geometrically reduced fibers of $Y$.
\end{proof}

\subsection{Existence of limits}

We assume that $S=\spec A$ is the spectrum of a discrete valuation
ring.  Let $Y/S$ be an arbitrary (not necessarily proper)
Deligne-Mumford stack of finite type, and let $V:=U\times_Y
U\rightrightarrows U\to Y$ be an \'etale presentation.  Let
$X_{\eta}\to Y_{\eta}$ be a finite morphism such that $X_{\eta}$ is
geometrically reduced.  Let $\ms R(X_{\eta})$ be the integral closure
of $\ms O_Y$ in $\ms O_{X_{\eta}}$, and let $\widetilde Y=\rspec\ms
R\to Y$.

\begin{lems}
  The formation of $\ms R$ commutes with \'etale base change $Y'\to
  Y$.
\end{lems}
\begin{proof}
  This follows immediately from 6.14.1 of \cite{ega4-2}.  Since this
  has a rather involved proof, we also offer a simpler alternative
  here.  It is a tautology that the formation of $\ms R$ commutes with
  Zariski base change.  Thus, it is easy to see that it suffices to
  prove this when $Y=\spec B$ is a local scheme and $Y'=\spec
  B[x]/(f(x))[1/f'(x)]$ is a basic \'etale morphism (so $f(x)$ is
  monic of some degree $n$).  Let $X=\spec C$.  Note that $B'=\spec
  B[x]/(f(x))$ is a finite free $B$-module.  Calculations due to Tate
  (which may be found in \S VII.1 of \cite{raynaud-hensel}) show that
  for any $y\in B'$, $f'(x)y=\sum_{i=0}^{n-1}\Tr_{B'/B}(b_iy)x^i$,
  with $b_1,\ldots,b_{n-1}$ certain elements of $B'$.  Let
  $C'=B'\tensor_B C$; since $B'/B$ is \'etale away from $Z(f'(x))$,
  $C'$ is reduced after inverting $f'(x)$.  If $z\in C'$ is any
  element integral over $B'$, then $\Tr_{C'/C}z$ is integral over $B$
  (\S 5.1, Prop.\ 17ff of \cite{bourbaki-comm-alg}).  Applying the
  formula, we see that if $y\in C'$ is integral over $B'$, then
  $f'(x)y\in B'$ (as it is a polynomial in $x$ with coefficients in
  $B$).  On the other hand, if $w$ is integral over $B'[1/f']$, then
  there is a multiple $(f')^sw$ which is integral over $B'$.  The
  result follows.
\end{proof}

\begin{lems}\label{L:branchvariety is normalization}
  If $X\to Y$ is a branchvariety, then the natural homomorphism
  $\rspec \ms R\to X$ is an isomorphism.
\end{lems}
\begin{proof}
  This is essentially Lemma 2.1 of \cite{alexeev-knutson}.  (While
  they work with graded rings, their proof carries over verbatim.)
\end{proof}

\begin{lems}\label{L:isoms}
  Suppose $X\to Y$ and $X'\to Y$ are branchvarieties.  The natural
  restriction map
$$\rho:\isom_Y(X,X')\to\isom_{Y_{\eta}}(X_{\eta},X'_{\eta})$$
is an isomorphism.
\end{lems}
\begin{proof}
  Given a generic isomorphism $\phi:X_{\eta}\simto X'_{\eta}$, there
  is an induced isomorphism $\rspec\ms R(X_{\eta})\to\rspec\ms
  R(X'_{\eta})$.  By \ref{L:branchvariety is normalization}, this
  gives rise to an isomorphism $X\to X'$.  Thus, $\rho$ is surjective.
  On the other hand, $\ms R(X_{\eta})\subset\iota_{\ast}\ms
  O_{X_{\eta}}$, where $\iota:X_{\eta}\to X$ is the inclusion of the
  generic fiber.  This immediately implies that $\rho$ is injective.
\end{proof}

\begin{cor}
  The stack of branchvarieties is separated with finite diagonal.
\end{cor}
\begin{proof}
  Lemma \ref{L:isoms} is precisely the valuative criterion of
  properness for the diagonal (when $X=X'$).  Quasi-finiteness of the
  diagonal follows immediately from the fact that the automorphism
  group of a finite reduced algebra over a field is finite.
\end{proof}

\begin{props}\label{P:descent of limit}
  There is a branchvariety $X\to Y$ extending $X_{\eta}$ if and only
  if there is a branchvariety $\mc X\to U$ extending $X_{\eta}\times_Y
  U$.
\end{props}
\begin{proof}
  The non-trivial part of the proposition is deducing the existence of
  $X$ from the existence of $\mc X$.  Consider the two pullbacks $\mc
  X_1:=p_1^{\ast}\mc X\to V$ and $\mc X_2:=p_2^{\ast}\mc X\to V$.
  Since $\mc X_{\eta}$ descends to $Y$, there is an isomorphism $\mc
  (X_1)_{\eta}\simto\mc (X_2)_{\eta}$ with trivial coboundary on
  $U\times_Y U\times_Y U$.  By \ref{L:isoms}, this descent datum
  extends to a descent datum $\mc X_1\simto\mc X_2$, yielding a finite
  $Y$-space $X\to Y$ and an isomorphism $X_U\simto\mc X$.  Since $U\to
  Y$ is \'etale, the fact that $\mc X\to U$ is a branchvariety
  immediately implies that $X\to Y$ is a branchvariety.
\end{proof}

\begin{lems}\label{L:affine limit}
  Suppose $Y=\spec R$ is an affine scheme of finite presentation over
  $S$.  Given a generic branchvariety $X_{\eta}\to Y_{\eta}$, there is
  a totally ramified extension $A\subset A'$ such that
  $X_{\eta}\tensor A'$ extends to a branchvariety over all of
  $Y\tensor A'$.
\end{lems}
\begin{proof}
  One way to prove this is to use Theorem 2.5 of
  \cite{alexeev-knutson}!  Choose an affine embedding
  $Y\subset\A^N_S$, and let $\widebar Y\subset\P^N_S$ be the
  projective closure of $Y$.  Normalizing $\widebar Y_{\eta}$ in $\ms
  O_{X_{\eta}}$ yields a generic branchvariety $\widebar
  X_{\eta}\to\widebar Y_{\eta}$ whose restriction to $Y_{\eta}$ is
  isomorphic to $X_{\eta}$.  By the result of Alexeev and Knutson
  cited in the first sentence, there is a totally ramified extension
  $A\subset A'$ and a branchvariety $\widebar X\to\widebar Y\tensor
  A'$ extending $\widebar X_{\eta}\tensor A'$.  Restricting $\widebar
  X$ to $Y\tensor A'$ yields the desired branchvariety $X\to Y\tensor
  A'$.
\end{proof}

\begin{proof}[Proof of \ref{T:limits}]
  We may assume that $S=\spec A$ is the spectrum of a discrete
  valuation ring.  Since $Y$ is of finite presentation over $S$, it is
  quasi-compact, so we may choose an \'etale cover $U\to Y$ with $U$
  affine.  Applying \ref{L:affine limit} and \ref{P:descent of limit}
  completes the proof.
\end{proof}

\subsection{Applications to branchvarieties on polarized
  orbifolds}\label{S:orbifolds}

The most potentially interesting application of these results is to
the study of branchvarieties on orbifolds with projective coarse
moduli spaces.  In this case, there will be enough numerical
invariants to again produce proper stacks.  More generally, there are
two (related) situations under which one can prove that the components
are quasi-compact.

\begin{thms}\label{T:cover means good}
  Suppose $S$ is Noetherian.  If there is a proper flat surjection
  $\mu:\Xi\to Y$ with $\Xi$ quasi-projective over $S$ then the
  connected components of $\br_{Y/S}$ are proper over $S$.
\end{thms}
\begin{proof}
  Note that the morphism sending $f:X\to Y$ to the coherent sheaf
  $f_{\ast}\ms O_X$ gives a finite-type morphism of stacks
  $\br_{Y/S}\to\coh_{Y/S}$.  Moreover, it sends a connected component
  into a connected component.  Furthermore, the pullback $\mu^{\ast}$
  gives a finite type homomorphism of algebraic stacks
  $\coh_{Y/S}\to\coh_{\Xi/S}$.  It thus suffices to show that the
  connected components of the stack $\coh_{\Xi/S}$ are of finite type
  over $S$.  But this follows from Lemma 1.3 and Lemma 1.6 of
  \cite{alexeev-knutson}.  Indeed, since $S$ is quasi-compact, there
  is an $S$-very ample invertible sheaf on $\Xi$ coming from a global
  immersion $\Xi\inj\P^N_S$.  Thus, it suffices to prove the result
  assuming $\Xi=\P^N_S$.  On any connected component of
  $\coh_{\P^N/S}$, the Hilbert polynomial and degree sequence are
  constant in fibers.  By Kleiman's theorem, such sheaves have a
  bounded Castelnuovo-Mumford regularity, so they all appear as
  quotients of a fixed sheaf of the form $\ms O(-m)^M$ on $\P^N_S$.
  The rest follows by classical results of Grothendieck.
\end{proof}

\begin{remarksub}
  It is tempting to believe that if $Y$ is an algebraic space of
  finite presentation over $S$ and there is a proper flat map $\Xi\to
  Y$ with $\Xi$ a quasi-projective scheme, then $Y$ is in fact a
  quasi-projective scheme.  When $S$ is the spectrum of a field, this
  is true.  The proof proceeds by (locally on $S$) slicing $\Xi$ to
  produce a finite flat such morphism (as in \cite{vistoli-kresch}).
  Taking the norm of an ample invertible sheaf on the cover then
  produces an ample invertible sheaf on $Y$ (cf.\ \S 6.6 of
  \cite{ega2}).  When $S$ is larger than a point, we can at least see
  that if $Y$ admits a finite flat cover by a quasi-projective scheme
  then $Y$ is quasi-projective.
\end{remarksub}

\begin{thms}\label{T:generating sheaf land}
  Suppose $S$ is Noetherian and $Y$ is a tame quotient stack with
  quasi-projective coarse moduli space.  The connected components of
  the stack $\br_{Y/S}$ are proper over $S$.
\end{thms}
\begin{proof}
  We can clearly assume that $S$ is connected.  If $S$ is the spectrum
  of a field, then this follows from the fact \cite{vistoli-kresch}
  that $Y$ admits a finite flat cover by a projective scheme, combined
  with \ref{T:cover means good}.  In general we do not know that $Y$
  admits a nice cover, but we can give an alternative proof.

  Let $E$ be a generating sheaf for $Y$, as defined in
  \cite{olsson-starr}.  Consider the morphism
  $\gamma_E:\coh_{Y/S}\to\coh_{\widebar Y/S}$ defined by sending $F$
  to $\pi_{\ast}\shom(E,F)$ (that $\gamma_E(F)$ is flat over the base
  follows from the tameness of $Y$).  By definition, there is a
  surjection $\pi^{\ast}\gamma_E(F)\tensor E\surj F$.  A connected
  component $\Gamma\subset\coh_{Y/S}$ maps into a connected component
  of $\coh_{\widebar Y/S}$.  By Kleiman's theorem, the connected
  components of $\coh_{\widebar Y/S}$ are proper.  More precisely,
  there exist $m$ and $M$ such that for all $F\in\Gamma$ there is a
  surjection $\ms O(-m)^M\surj\pi_{\ast}\shom(E,F)$.  Combining this
  with the definition of a generating sheaf, this yields a surjection
  $E(-m)^M\surj F$.  Moreover, since $\Gamma$ is connected, for any
  $T\to S$ and any lift $T\to\Gamma$, the function
  $P_{F_t}:K^0(Y)\to\Z$ sending a locally free sheaf $G$ on $Y$ to
  $\chi(\shom(G_t,F_t))$ is locally constant as $t$ varies (Lemma 4.3
  of \cite{olsson-starr}).  Thus, $\Gamma$ is the image of a closed
  subspace of the space $Q(E(-m)^M,P)$ of quotients of $E$ with
  Hilbert polynomial $P$.  By 4.5 of \cite{olsson-starr}, this space
  of quotients is quasi-projective over $S$, hence is quasi-compact.
  (The hypothesis in \cite{olsson-starr} that $S$ be an affine scheme
  is unnecessary for $Q(P)$ to be a quasi-compact algebraic space, but
  the reader uncomfortable with this bald assertion may choose to only
  regard this corollary as true under the additional hypothesis that
  $S$ is an affine scheme.)  It follows that $\Gamma$ is
  quasi-compact, finishing the proof.
\end{proof}

Of course, describing the connected components (or collections
thereof) is quite a subtle task.  The best one could hope for is to
find numerical invariants which bound a collection of branchvarieties.
As the proofs above show, once one has chosen a polarized cover or a
generating sheaf and a polarization of the coarse moduli space, one
can use the resulting numerical invariants to bound substacks of
$\br_{Y/S}$.

If $S$ is the spectrum of a field and $Y$ is tame, smooth, and
separated with quasi-projective coarse moduli space, then it is known
\cite{vistoli-kresch} that $Y$ is a quotient stack (and that it admits
a finite flat cover by a quasi-projective scheme).  Thus, \ref{T:cover
  means good} and \ref{T:generating sheaf land} both apply to show
that the components of $\br_{Y/S}$ are proper.  In this direction, it
might be interesting to consider the scheme of branchvarieties of a
weighted projective space.  The corresponding Hilbert scheme plays an
important role in recent work of Abramovich and Hassett on the moduli
of stable varieties.


\begin{thebibliography}{1}

\bibitem{alexeev-knutson}
Valery Alexeev and Allen Knutson.
\newblock Complete moduli spaces of branchvarieties, 2006.
\newblock Preprint.

\bibitem{bourbaki-comm-alg}
Nicolas Bourbaki.
\newblock {\em Commutative algebra. {C}hapters 1--7}.
\newblock Elements of Mathematics (Berlin). Springer-Verlag, Berlin, 1998.
\newblock Translated from the French, Reprint of the 1989 English translation.

\bibitem{ega2}
A.~Grothendieck.
\newblock \'{E}l\'ements de g\'eom\'etrie alg\'ebrique. {II}. \'{E}tude
  globale \'el\'ementaire de quelques classes de morphismes.
\newblock {\em Inst. Hautes \'Etudes Sci. Publ. Math.}, (8), 1961.

\bibitem{ega3-2}
A.~Grothendieck.
\newblock \'{E}l\'ements de g\'eom\'etrie alg\'ebrique. {III}. \'{E}tude
  cohomologique des faisceaux coh\'erents. {II}.
\newblock {\em Inst. Hautes \'Etudes Sci. Publ. Math.}, (17):223, 1963.

\bibitem{ega4-2}
A.~Grothendieck.
\newblock \'{E}l\'ements de g\'eom\'etrie alg\'ebrique. {IV}. \'{E}tude locale
  des sch\'emas et des morphismes de sch\'emas. {II}.
\newblock {\em Inst. Hautes \'Etudes Sci. Publ. Math.}, (24):231, 1965.

\bibitem{ega4-3}
A.~Grothendieck.
\newblock \'{E}l\'ements de g\'eom\'etrie alg\'ebrique. {IV}. \'{E}tude locale
  des sch\'emas et des morphismes de sch\'emas. {III}.
\newblock {\em Inst. Hautes \'Etudes Sci. Publ. Math.}, (28):255, 1966.

\bibitem{vistoli-kresch}
Andrew Kresch and Angelo Vistoli.
\newblock On covering of {D}eligne-{M}umford stacks and surjectivity of the
  {B}rauer map.
\newblock {\em Bull. London Math. Soc.}, 36(2):188--192, 2004.

\bibitem{lmb}
G\'erard Laumon and Laurent Moret-Bailly.
\newblock Champs alg\'ebriques.
\newblock Springer-Verlan, Berlin, 2000.

\bibitem{olsson-sheaves}
Martin Olsson.
\newblock Sheaves on {A}rtin stacks.
\newblock Preprint, 2005.

\bibitem{olsson-homs}
Martin Olsson.
\newblock Hom--stacks and restriction of scalars.
\newblock To appear, {\em Duke Math.\ J.}.

\bibitem{olsson}
Martin Olsson.
\newblock On proper coverings of {A}rtin stacks.
\newblock {\em Adv. Math}, 198: 93--106, 2005.

\bibitem{olsson-starr}
Martin Olsson and Jason Starr.
\newblock Quot functors for {D}eligne-{M}umford stacks.
\newblock {\em Comm. Algebra}, 31(8):4069--4096, 2003.
\newblock Special issue in honor of Steven L. Kleiman.

\bibitem{raynaud-hensel}
Michel Raynaud.
\newblock {\em Anneaux locaux hens\'eliens}.
\newblock Lecture Notes in Mathematics, Vol. 169. Springer-Verlag, Berlin,
  1970.

\bibitem{starr}
Jason Starr.
\newblock Remarks on moduli spaces and Artin's axioms.
\newblock Preprint, 2006.

\end{thebibliography}
\end{document}